\documentclass[12pt]{article}
\usepackage[dvips]{graphicx}
\usepackage{graphicx, color}
\usepackage{tikz}
\usepackage{hyperref}
\usepackage{framed}
\usepackage{amsfonts,amsmath, amssymb, amstext, amsthm}

\DeclareGraphicsExtensions{ps,eps,pdf}
\theoremstyle{plain}
\newtheorem{theorem}{Theorem}

\newtheorem{lemma}{Lemma}
\theoremstyle{remark}

\theoremstyle{definition}
\newtheorem{definition}{Definition}

\def\R{{\mathbb R}}

\def\cF{{\cal F}}

\def\cM{{\cal M}}

\def\endproof{\hfill $\square$\bigskip}
\def\newchi{\raise2pt\hbox{$\chi$}}
\def\newdag{\raise4pt\hbox{\dag}}

\begin{document}

\title{Image Processing Variations with Analytic Kernels
\thanks{We thank Inwon Kim and Peter Petersen for helpful conversations.  The first named author thanks the CRM Barcelona for support during its
 Research Program on Harmonic Analysis,
Geometric Measure Theory and Quasiconformal Mappings, where this
 paper was presented and where part of its research was undertaken.}  \thanks{This work was supported in part by the National  Science Foundation (Grant NSF DMS 071495 and 0809270) and from the Office of Naval Research (Grant ONR N00014-09-1-0108).}}
 \author{John B. Garnett \thanks{Department of Mathematics, University of California, Los  Angeles, 405 Hilgard Ave., Los Angeles, CA 90095-1555, USA ({\tt jbg@math.ucla.edu}).}  \and Triet M. Le \thanks{Department of Mathematics, University of Pennsylvania, Philadelphia, PA 19104, USA ({\tt trietle@math.upenn.edu}).}
\and Luminita A. Vese \thanks{Department of Mathematics, University of California, Los Angeles, 405 Hilgard Ave., Los Angeles, CA 90095-1555, USA ({\tt lvese@math.ucla.edu}).}}
 
 \date{}  

\maketitle
\thispagestyle{plain}

\begin{abstract}
Let $f\in L^1(\R^d)$ be real. The Rudin-Osher-Fatemi model is to minimize 
$\|u\|_{\dot{BV}}+\lambda\|f-u\|_{L^2}^2$,
in which  one thinks of $f$ as a given image, $\lambda > 0$ as a ``tuning parameter",  $u$ as an optimal ``cartoon''  approximation to $f$, and $f-u$ as ``noise'' or ``texture''. 
Here we study variations of the R-O-F model having the form  
 $\inf_u\{\|u\|_{\dot{BV}}+\lambda \|K*(f-u)\|_{L^p}^q\}$ where $K$ is a real analytic kernel such as a Gaussian.  
For these functionals we characterize the minimizers  $u$  and establish several  of their properties, including especially
their smoothness properties. In particular we prove that on any open set 
on which $u \in W^{1,1}$ and $\nabla u \neq 0$ almost every level set $\{u =c\}$ is a real analytic surface.  
We also prove 
that if $f$ and $K$ are radial functions then every minimizer $u$ is a radial step function.

\end{abstract}


\section{Introduction}


 Several $BV$ variational models have been proposed as image decomposition models (see Section 2 for the definition of $BV$).
First,  Rudin-Osher-Fatemi \cite{ROF} proposed the minimization 
\begin{equation}\label{rof}
\inf_{u\in BV} \left\{\|u\|_{\dot{BV}}  + \lambda \|f-u\|^2_{L^2}\right\}.
\end{equation} 
In (\ref{rof}), $f\in L^1(\R^d)$ is a real function and one thinks of $u$ as the  ``cartoon'' component of $f$ and $f-u$ as the  ``noise+texture'' component of $f$. By the strict convexity of the  
functional $\|f-u\|^2_{L^2}$, problem  (\ref{rof}) has  a unique minimizer $u$.  However, one limitation of model (\ref{rof}) is illustrated by the following example from \cite{Meyer} and \cite{ChanStrong}:  
if $d =2$ and $f=\alpha\chi_D$ where $D$ a disk centered at the origin and of radius $R$, then $u=(\alpha-(\lambda R)^{-1})\chi_D$ and $v=f-u=(\lambda R)^{-1}\chi_D$ if $\lambda R\geq 1/{\alpha}$, but $u =0$ 
 if $\lambda R\leq 1/{\alpha}$. Thus $u\neq f$ can occur even though $f\in BV$ is already a cartoon without texture or noise (note that $f$ and $u$ still have the same set of discontinuity). 
To overcome this limitation and also to attempt to separate noise from texture, many authors 
have introduced alternate forms  of (\ref{rof})
by replacing $\|f-u\|^2_{L_2}$ by other expressions. We mention the book \cite{Meyer} and the papers  
  \cite{TNV1}, \cite{TNV2},  \cite{TadmorAthavale},   
 \cite{ChanEsedoglu}, 
 \cite{Allard1,Allard2,Allard3},    
 \cite{VO}, 
 \cite{OSV}, \cite{LinhLieu}, \cite{LeVese},
 \cite{Aujol1,Aujol2}, 
 \cite{GLMV},  \cite{AujolChambolle}),  \cite{GJLV},
 \cite{aujolHilbert}.
Among these, the papers of Chan and Esedoglu \cite{ChanEsedoglu} and Allard  
\cite{Allard1,Allard2,Allard3} are closest to the present work.

Chan and Esedoglu \cite{ChanEsedoglu}  considered  the minimization
$$
\inf_{u\in BV}\Big\{\|u\|_{\dot{BV}}+\lambda\int|f-u|dx\Big\}
$$
 (see also Alliney \cite{Alliney} for the one-dimensional discrete case).
For this problem  minimizers always exist but they may not be unique. 
For the example $d=2$ and $f=\chi_{B(0,R)},$
\cite{ChanEsedoglu} gives  
$u=f$ if $R> { {2}\over {\lambda} } $ and $u=0$ if $R< { {2}\over {\lambda} } $.
W. Allard \cite{Allard1,Allard2,Allard3} analyzed extremals for the problem 
$$\inf_{u\in BV}\Big\{\|u\|_{\dot{BV}}+\lambda\int \gamma(u-f)dx\Big\}$$
where
$\gamma(0)=0$, $\gamma \geq 0$, and $\gamma$ is locally Lipschitz. 
Then minimizers $u$ exist although they may not be unique.  Moreover, the minimizers $u$
satisfy the smoothness condition
$$\partial^*(\{u>t\})\in C^{1+\alpha},\ \ \ \alpha\in (0,1)$$
where $\partial^*$ denotes ``measure theoretic boundary". 
 Allard also gave mean curvature estimates on $\partial^*(\{u>t\}).$

In this paper we study a cartoon+texture decomposition model  defined with a positive, real analytic  convolution kernel $K$:  
\begin{equation}\label{(1.2)}
\inf_{u\in BV} \Big\{\|u\|_{\dot{BV}}  + \lambda \|K*(f-u)\|^q_{L^p}\Big\}
\end{equation}  where  $1 \leq p, q < \infty.$ 
We choose the kernel $K$ in \eqref{(1.2)} so that the Fourier transform
 $\widehat{K}(\xi)$ decays
rapidly as $|\xi|\rightarrow \infty$. The motivation is that we expect $v=f-u$ to be oscillatory, so that $\widehat{v}(\xi)$ is large when $|\xi|$ is large. Thus, $\widehat{K}\cdot\widehat{v}= \widehat{(K * v)}$ dampens high frequencies of $v$, which suggests that $\|K*v\|_{L^p}^q$ is small for oscillatory $v$. 
We also want the cartoon component $u$ 
to be very simple, for example, to be piecewise constant or to have real analytic level sets, and 
for that reason we choose $K$ to be real analytic.  Examples of such $K$ are the Gaussian kernel where $\widehat{K}(\xi) = e^{-\pi t|\xi|^2}$ or the Poisson kernel where $\widehat{K}(\xi) = e^{-\pi t|\xi|}$, for some $t>0$.


By comparison  \cite{ChanEsedoglu} takes $p = q =1$ and $K = {\rm identity}$ and
our choices of $K$ yield more precise results about the minimizers for (\ref{(1.2)}).  In comparison with Allard's paper \cite{Allard1} we note that for many
 choices
of the kernel $K$ our functional $||K*(f-u)||^q_{L^p}$ is {\em admissible} in the sense of \cite{Allard1}
so that 
  the regularity  results from section 1.5 of that paper  hold for the minimizers $u$ of (\ref{(1.2)}).
However,  because of the analyticity of $K$ our minimizers have greater  smoothness than those from 
 \cite{Allard1}.   Moreover the functional in (\ref{(1.2)}) is not {\em local} in the sense of \cite{Allard1}, so that the conclusions 
of section 1.6 of \cite{Allard1} need not hold for the  minimizers of (\ref{(1.2)}). 



\section{The Variational Problems}

To begin we recall the definition of $BV = BV(\R^d).$ 

\begin{definition}
Let $u \in L^1_{\rm {loc}}(\R^d)$ be real.   
  We say $u \in BV$ if 
$$\sup \Bigl\{\int u {\rm  {div}}\varphi  dx: \varphi \in C^1_0(\R^d), \sup|\varphi(x)| \leq 1\Bigr\} 
=\|u\|_{\dot{BV}} < \infty.$$ 
\end{definition}

\noindent If $u \in BV$ there is an $\R^d$-valued measure $\vec \mu$ such that
 ${{\partial u} \over {\partial x_j}} = (\vec \mu)_j$ as distributions and we write
$$Du = \vec \mu.$$
The  vector
measure $\mu$ has a {\it {polar decomposition}} 
$$\vec \mu = \vec \rho \mu$$ 
where $\mu$ is a finite positive Borel measure and $\vec \rho:\R^d\to S^{d-1}$ is a Borel function, and 
$$\|u\|_{\dot{BV}} = \int d\mu. $$

\noindent(see for example Evans-Gariepy \cite{EvansGariepy}).

We assume  $K$ is a  positive,
even, bounded and real analytic kernel on $\R^d$ such that 
$\int Kdx =1$ and such that $K*u$ determines $u$ (i.e. the map $L^p \ni u \to K*u$ is injective).
For example we may take $K$ to be a Gaussian or a Poisson kernel.   
We fix  $\lambda > 0$,   $1 \leq p < \infty$ and $1 \leq q < \infty.$  
For real $f(x) \in L^1$ we
 consider the extremal problem:
\begin{equation}
m_{p,q,\lambda}  = \inf\{\|u\|_{\dot{BV}} + \cF_{p,q,\lambda}(f-u): u \in BV\} 
\label{eq(2.1)}
\end{equation}
where
\begin{equation}
\cF_{p,q,\lambda}(h)  = \lambda \|K * h\|^q_{L^p}. 
\label{eq(2.2)}
\end{equation}
 Since $BV \subset L^{{d} \over {d-1}}$ and $K \in L^{\infty}$, a weak-star compactness argument shows that (\ref{eq(2.1)}) has at least one minimizer $u$.  Our objective is to describe, given $f$,  the set $\cM_{p,q,\lambda}(f)$ of minimizers $u$ of (\ref{eq(2.1)}). 



\subsection{Convexity}

Since the functional in (\ref{eq(2.1)})  is convex, the set of minimizers 
  $\cM_{p,q,\lambda}(f)$ is a convex subset of $BV$.   If $p > 1$ or 
if $q > 1$, then the functional
(\ref{eq(2.2)}) is strictly convex and the problem (\ref{eq(2.1)}) has a unique minimizer because $K*u$ determines $u$. 
When $p = q =1$ minimizers may not be unique, but they satisfy the relations  
given in (\ref{eq(2.3)}) and (\ref{eq(2.4)}) below. 

\begin{lemma} 
\label{Lemma2.1}
Let $p = q =1$  and assume  $u_1 \in \cM_{p,q,\lambda}(f)$ and 
$u_2 \in \cM_{p,q,\lambda}(f)$.
For $j =1,2$ write
$$Du_j = \vec \mu_j = \vec \rho_j \mu_j$$ 
with $|\vec \rho_j| =1$ and $\mu_j \geq 0$ and 
write ${{d\vec \mu_j} \over {d\mu_k}}$
for the Radon-Nikodym derivative of (the absolutely continuous part of) $\vec \mu_j$ 
with respect to $\mu_k.$  
 Then
\begin{equation}
{{K*(f - u_1)} \over {|K*(f-u_1)|}} ={{K*(f-u_2)} \over {|K * (f-u_2)|}} 
~~~almost ~ everywhere 
\label{eq(2.3)}
\end{equation}
on $\{|K * (f -u_j)| > 0\}, j = 1, 2;$
\noindent  and  
\begin{equation}
\vec \rho_k \cdot {{d\vec \mu_j} \over {d\mu_k}} = \Big|{{d\vec \mu_j} \over {d\mu_k}}\Big|, ~j \neq k. 
\label{eq(2.4)}
\end{equation}
\end{lemma}

\noindent {\bf Proof:} Since $\cM_{p,q,\lambda}(f)$ is a convex subset of $BV$, ${{u_1 + u_2} \over {2}}$ is also a minimizer. This implies
\begin{equation}\label{lem_eq_1}
\begin{split}
\left\|\frac{u_1+u_2}{2}\right\|_{\dot{BV}} +\lambda \left\|K*\left(f-\frac{u_1+u_2}{2}\right)\right\|_1 &= \frac{1}{2}\bigl(\|u_1\|_{\dot{BV}}
 + \|u_2\|_{\dot{BV}}
\bigr) \\
&+ \frac{\lambda}{2}\bigl(\|K * (f -u_1)\|_1 + \|K * (f -u_2)\|_1\bigr).
\end{split}
\end{equation}
On the other hand, using the convexity of $\|\cdot\|_{\dot{BV}}
$ and $\|\cdot\|_{L^1}$ we have
\begin{equation}\label{lem_eq_2}
\left\|\frac{u_1+u_2}{2}\right\|_{\dot{BV}}
 \le \frac{1}{2}\bigl(\|u_1\|_{\dot{BV}}
 + \|u_2\|_{\dot{BV}}
\bigr)
\end{equation}
\noindent and 
\begin{equation}\label{lem_eq_2a} 
\left\|K*\left(f-\frac{u_1+u_2}{2}\right)\right\|_1\le \frac{1}{2}\bigl(\|K * (f -u_1)\|_1 + \|K * (f -u_2)\|_1\bigr)
\end{equation}
Combining \eqref{lem_eq_1}, \eqref{lem_eq_2}, and \eqref{lem_eq_2a} we obtain
the equality

\begin{equation}\label{lem_eq_2.5}
\Big|\Big|K*(f - {{u_1 + u_2} \over {2}})\Big|\Big|_1 = {{1} \over {2}} 
\bigl(\|K * (f -u_1)\|_1 + \|K*(f -u_2)\|_1\bigr),
\end{equation}
which implies (\ref{eq(2.3)}). We also obtain 
\begin{equation}\label{lem_eq_3}
\left\|u_1+u_2\right\|_{\dot{BV}}
 = \|u_1\|_{\dot{BV}}
 +\|u_2\|_{\dot{BV}}
\end{equation}
and for  $k\neq j$ equation \eqref{lem_eq_3} implies
$$\int \Big|\vec{\rho}_k + {{d\vec \mu_j} \over {d\mu_k}}\Big| d\mu_k = \int d\mu_k + \int\Big|{{d \vec \mu_j}
\over {d\mu_k}}\Big| d\mu_k,$$
which yields (\ref{eq(2.4)}). 
\endproof

\subsection{Properties of $u \in \cM_{p,q,\lambda}(f)$}

\begin{lemma} 
\label{Lemma2.2}
Let  $f\in L^1$,  let $u \in BV$ be a minimizer of (\ref{eq(2.1)}) with $\|u-f\|_1 \neq 0$ and write    
$$
Du = \vec{\mu} = \vec{\rho}\cdot \mu.
$$
 Then whenever $h \in BV$ is real,  $Dh = \vec{\nu}$ and 
 $ \vec \nu = {{d\vec \nu} \over {d\mu}} \mu + \vec \nu_s$
 is the Lebesgue decomposition of $\vec \nu$ with respect to $\mu$ 
(so that $\vec \nu_s$ is singular to $\mu$), we have  
\begin{equation}
\left| \int \vec{\rho}  \cdot  {{d \vec \nu} \over {d\mu}} d\mu 
- \lambda  \int h (K * J_{p,q}) dx\right| \leq \|\vec \nu_s\|,  
\label{eq(3.1)}
\end{equation}
where 
\begin{equation}
 J_{p,q} = q{{F |F|^{p-2}} \over {\|F\|_p^{p-q}}},
 ~~  F = K *(f-u)
\label{eq(3.2)}
\end{equation}
and $\|\vec \nu_s\|$ denotes the
 norm of the vector measure $\vec \nu_s$.
 Conversely, if $u \in BV$,  $\|u -f\|_1 \neq 0$ and if (\ref{eq(3.1)}) and (\ref{eq(3.2)}) hold for all $h$,
 then $u \in \cM_{p,q,\lambda}(f).$ 
\end{lemma}
\bigskip
\noindent Note that because $\|u - f\|_1 \neq 0$ and  $K *(f-u)$
is real analytic and bounded, $J_{p,q}$ is defined almost everywhere,
and that by Lemma \ref{Lemma2.1}, $J_{p,q}$ is independent of $u \in \cM_{p,q,\lambda}$
in the case $p = q =1$. 

\bigskip
\noindent {\bf {Proof:}} Let $|\epsilon|$ be sufficiently small.  Since
$u$ is extremal, we have
\begin{equation}\label{lem2_eq1}
\|u + \epsilon h\|_{\dot{BV}
} - \|u\|_{\dot{BV}} + \cF_{p,q,\lambda}(f - u - \epsilon h)-
\cF_{p,q,\lambda}(f -u) \geq 0.
\end{equation}
On the other hand we have
$$
\left|\vec{\rho}  + \epsilon  {{d\vec{\nu}} \over {d\mu}}\right| = \left(1 + 2\epsilon \vec{\rho}\cdot \frac{d\vec{\nu}}{d\mu} + \epsilon^2\left\|\frac{d\vec{\nu}}{d\mu} \right\|^2\right)^{1/2} =  1 + \epsilon \vec{\rho}\cdot \frac{d\vec{\nu}}{d\mu} + o(|\epsilon|),
$$
where in the last equality, we use the estimate $(1+\alpha)^{1/2} = 1 + \frac{\alpha}{2} + o(|\alpha|)$. This implies,
\begin{equation*}
\begin{split}
\|u + \epsilon h\|_{\dot{BV}} - \|u\|_{\dot{BV}} =|\epsilon| \|\vec{\nu}_s\|
+ \int \left(\left|\vec{\rho}  + \epsilon 
 {{d\vec{\nu}} \over {d\mu}}\right| - 1\right) d\mu  =
|\epsilon| \|\vec{\nu}_s\| + \epsilon \int \vec{\rho} \cdot  
{{d\vec{\nu}} \over {d\mu}} d\mu + o(|\epsilon|).
\end{split}
\end{equation*}
Moreover $K * (f -u)$ is bounded and non-zero almost everywhere, since $K$ is real analytic.  Hence we also have 
\begin{equation*}
\begin{split}
\cF_{p,q,\lambda}(f -u -\epsilon h) - \cF_{p,q,\lambda}(f-u) &=
 -\lambda \epsilon \int(K * h) J_{p,q} dx + o(|\epsilon|) \\
&= - \lambda \epsilon \int h (K*J_{p,q}) dx + o(|\epsilon|) 
\end{split}
\end{equation*}
since $K(-x) = K(x)$. Thus by \eqref{lem2_eq1}, we have
$$
-\epsilon \left[ \int \vec{\rho} \cdot  
{{d\vec{\nu}} \over {d\mu}} d\mu - \lambda  \int h (K*J_{p,q}) dx\right] \le |\epsilon|\|\vec{\nu}_s\| + o(|\epsilon|). 
$$
Taking $\pm \epsilon$ and noting that the right side of the above inequality does not depend on the sign of $\epsilon$, we see that  (\ref{eq(3.1)}) holds. 
The converse statement holds because the functional (\ref{eq(2.2)}) is convex. 
\endproof

Lemma  \ref{Lemma2.2} does not hold  for the Chan-Esedoglu \cite{ChanEsedoglu}
functional because in that case one can have $f-u =0$ on a set of positive measure, and this yields 
the additional term $\int_{\{|f -u| = 0\}} |h| dx$ on the right side of (\ref{eq(3.1)}).
 
\medskip

Later we will  need the following alternate characterization of minimizers, due to 
Meyer \cite{Meyer}  in the case of the Rudin-Osher-Fatemi model. 
Define 
$$
\|v\|_* = \inf\left\{\left\| |u|\right\|_{\infty}: v = \sum_{j=1}^d{{\partial u_j} \over {\partial x_j}}, |u|^2 = \sum_{i=1}^d |u_j|^2\right\}$$
\noindent so that $\|v\|_*$ is (isometrically) the norm of the dual of $W^{1,1}\subset BV$ 
when $W^{1,1}$ is given the norm of $BV$.  
By the weak-star density of $W^{1,1}$ in $BV$, 
\begin{equation}
\left|\int hv dx\right| \leq \|h\|_{\dot{BV}}
 \|v\|_* 
\label{eq(3.4)}
\end{equation}
whenever $v \in L^2$.
The lemma characterizes minimizers in  terms of $\|\cdot\|_*$. 

\begin{lemma} 
\label{Lemma2.3}
Let  $u  \in BV$ such that $u \neq f$, and let $J_{p,q}$ be defined as in Lemma 2.  
Then $u$ is a minimizer for the problem (\ref{eq(2.1)}) if and only if
\begin{equation}
\|K * J_{p,q}\|_* = {{1} \over {\lambda}}
\label{eq(3.5)}
\end{equation}
and 
\begin{equation}
\int u (K * J_{p,q}) dx = {{1} \over {\lambda}} \|u\|_{\dot{BV}}. 
\label{eq(3.6)}
\end{equation}

\end{lemma} 

\noindent {\bf {Proof:}}  The short proof is the same as
in \cite{Meyer}, but we include it for the reader's convenience.  Let  $u$ is a minimizer for (\ref{eq(2.1)}). Then for any $h \in W^{1,1}$, (\ref{eq(3.1)}) yields 
$$\left|\int h (K * J_{p,q}) dx\right| \leq {{\|h\|_{\dot{BV}}
} \over {\lambda}}$$
by the definition of $\vec \nu_s$. Hence by the definition of $\|\cdot\|_*$, 
$$\|K *J_{p,q}\|_* \leq {{1} \over {\lambda}}.$$
But setting $h =u$ in (\ref{eq(3.1)}) gives (\ref{eq(3.6)}), so that (\ref{eq(3.5)}) follows. 

Conversely, assume  $u \in BV$ satisfies (\ref{eq(3.5)}) and (\ref{eq(3.6)})  and note
 that $u$ determines $J_{p,q}$.  
Still following Meyer \cite{Meyer}, we let $h \in BV$ be real.  Then for small $\epsilon > 0$,
(\ref{eq(3.4)}), (\ref{eq(3.5)}) and (\ref{eq(3.6)}) give
\begin{eqnarray*}
\|u + \epsilon h\|_{\dot{BV}}
 + \lambda \|K * (f - u -\epsilon h)\|_1 
&\geq& \lambda \int (u + \epsilon h) (K * J_{p,q})dx
+ \lambda \|K * (f -u)\|_1 \\
&-&\epsilon \lambda \int h (K * J_{p,q}) dx + o(\epsilon)\\
&=&  \|u\|_{\dot{BV}} + \epsilon \lambda \int h (K * J_{p,q}) dx  \\
& -&\epsilon \lambda \int h (K *J_{p,q})dx + o(\epsilon) \\
&\geq& 0.
\end{eqnarray*} 
Therefore $u$ is a local minimizer for the functional (\ref{eq(2.1)}),
and by convexity that means $u$ is a global minimizer. 
\endproof
\medskip
 
\begin{lemma}
\label{Lemma2.4}
Assume  $f \in L^1, ~  u \in \cM_{p,q,\lambda}(f)$, and $\|u -f\|_1 \neq 0$. Let $U$ be an open set
on which  $Du = \vec \mu$ is absolutely continuous to Lebesgue measure and has
Radon-Nikodym derivative ${{d \vec \mu} \over {dx}} \neq 0$ almost everywhere.  
  Then as distributions on $U$ 
\begin{equation}
{\rm{div}}\left( {{\nabla u} \over {|\nabla u|}} \right)= -\lambda K *J_{p,q},
\label{eq(Lap)}
\end{equation}
and $u \in W^{1,1}(U).$ 
In particular, if $u \in C^2(U)$ then the level set $\{u =c\}$ is 
locally a 
$C^2$ surface having   
mean curvature 
$-\lambda K * J_{p,q}(x)$
at $x \in U.$
\end{lemma}
\noindent {\bf Proof:} Since $Du$ is absolutely continuous on $U$ we
have  $u \in W^{1,1}(U)$ and $\vec \mu = \nabla u dx$ there. 
 Let $h \in C^{\infty}$ have compact support contained in 
$U$. Then by the hypotheses, 
$\vec \nu = Dh = \nabla h dx$ is absolutely continuous to $Du$ so that by (\ref{eq(3.1)})
\begin{equation}
\label{eq(Lap1)}
\int_U \nabla h \cdot {{\nabla u} \over {|\nabla u|}} dx  = \lambda \int_U h (K * J_{p,q}) dx. 
\end{equation} 
This implies (\ref{eq(Lap)}).   Also,  if $u \in C^2(U)$ then 
(\ref{eq(Lap1)}) holds pointwise and gives the mean curvature of $\{u =c\}$ inside
$U$. 
\endproof

Known results on mean curvature equations can now be used to show that  almost every level set
$U \cap \{u =c\}$ is a real analytic surface, even without the assumption $u \in C^2(U)$.  Below we write  $\Lambda_{d-1}$ 
for $d-1$ dimensional Hausdorff measure.

\begin{theorem}
\label{Thm2.2}
Assume  $f \in L^1, ~  u \in \cM_{p,q,\lambda}(f)$, and $\|u -f\|_1 \neq 0$.
Let $U$ be an open set
on which  $Du = \vec \mu$ is absolutely continuous to Lebesgue measure and 
on which the  
Radon-Nikodym derivative ${{d \vec \mu} \over {dx}} \neq 0$ almost everywhere.
  Then for almost all $c \in \R$  and for $\Lambda_{d-1}$  almost every 
$x_0 \in U \cap \{u =c\}$  there exists a 
 $C^1$-hypersurface $S$ with continuous unit normal
$\vec n(x) = {{\nabla u} \over {|\nabla u|}}$ and a neighborhood 
$V$ of $x$ such that $\Lambda_{d-1}((V \cap \{u =c\}) \Delta S) = 0.$ 
After a rotation 
$S = \{x_d = \varphi(y): ~ y = (x_1, \dots, x_{d-1}) \in V_0\}$, where
$V_0 \subset {\R}^{d-1}$ is open,  $\varphi \in C^1(V_0)$, 
and $\vec n(y, \varphi(y)) = (1 + |\nabla \varphi|^2)^{-1/2} (\nabla \varphi, -1).$  Moreover,
as a distribution on $V_0$ 
\begin{equation}
\label{eq(Nab)}
{\rm {div}} \left(  {{\nabla \varphi} \over {(1 + |\nabla \varphi|^2)^{1/2}}}\right)
 = -\lambda K *J_{p,q}(y,\varphi(y))dy,
\end{equation}
and the function $\varphi$ and the surface $S$ are real analytic.

\end{theorem}

\noindent {\bf Proof:} That $S$ and $\varphi$ exist almost everywhere 
follows from standard properties of BV functions and the hypothesis
that $|\nabla u| > 0$ a. e. on $U$.  See the proof of Theorem 4 below and 
Chapter 5 of  \cite{EvansGariepy}. To prove (\ref{eq(Nab)}) we may 
assume $c =0$. Let
$h \in C^{\infty}_0(V_0)$, let $\newchi_{\epsilon}(t) = {{1} \over {\epsilon}}
\newchi({{t} \over {\epsilon}})$ where $\newchi(t) = \newchi(-t) \geq 0$
is $C^{\infty}(-1,1)$ and $\int \newchi dt =1,$ and define 
$$H_{\epsilon}(x) = \newchi_{\epsilon}(h(x_1, \dots, x_{d-1}) - x_d) h(x_1, \dots,
x_{d-1}).$$
Then by (\ref{eq(Lap)}),
\begin{eqnarray*}
& &\int\Biggl( \sum_{j=1}^{d-1}\bigl(\newchi_{\epsilon}(h(x_1, \dots, x_{d-1}) - x_d)
+ \newchi_{\epsilon}^{'}(h(x_1, \dots, x_{d-1}) - x_d)h(x_1, \dots, x_{d-1})\bigr)
\\
& &{{\partial h} \over {\partial x_j}} {{1} \over {|\nabla u|}} {{\partial u}
\over {\partial x_j}}\Biggr) -
\Biggl( \newchi_{\epsilon}^{'}(h(x_1, \dots, x_{d-1}) - x_d)h(x_1, \dots, x_{d-1})  {{1} \over {|\nabla u|}} {{\partial u}
\over {\partial x_d}}\Biggr) dx\\ &=& \lambda \int_{V_0}
 H_{\epsilon}(x) K*J_{p,q}(x) dx.  
\end{eqnarray*}
Now for almost every $c$ the right side of this equation tends
to $\lambda \int_{V_0} h (K * J_{p,q})(y) dy$
and, by the fine properties of BV functions 
in Chapter 5 of  \cite{EvansGariepy} or Chapter 3 of 
\cite{AmbrosioFuscoPallara}, the left side  tends to
$$\int_{V_0} \nabla h \cdot 
 {{\nabla \varphi} \over {(1 + |\nabla \varphi|^2)^{1/2}}}dy.$$
That proves (\ref{eq(Nab)}). 

To prove the real analyticity of $\varphi$, and hence of $S$, we invoke three theorems. First, since $\varphi \in C^1$, the results on mean curvature equations in 
 Section 7.7 of \cite{AmbrosioFuscoPallara} show that $\varphi \in W^{2,2} \cap C^{1 + \alpha}$
whenever $0 < \alpha < 1.$    
Next, since $\varphi \in W^{2,2}$ we can rewrite (\ref{eq(Nab)}) as 

\begin{equation}
\label{eq(Nab1)}
\sum_{j,k} {{\delta_{j,k} - \varphi_j \varphi_k} \over 
{(1 + |\nabla \varphi|^2)^{3/2}}} \varphi_{j,k} = \lambda
K*J_{p,q}(y,\varphi(y)).
\end{equation}
Indeed,  (\ref{eq(Nab1)}) is clear if $\varphi \in C^2$, and if we set $\varphi^{\epsilon}
= \newchi_{\epsilon} *\varphi \in C^2$ then in the norms of $C^{1 + \alpha}$ and  
$W^{2,2}$, $\varphi^{\epsilon} \to \varphi$ as $\epsilon \to 0.$  Hence for each $j$
$$\int_{V_0} h_j \sum_k  {{\delta_{j,k} - \varphi^{\epsilon}_j \varphi^{\epsilon}_k} \over 
{(1 + |\nabla \varphi^{\epsilon}|^2)^{3/2}}} \varphi^{\epsilon}_{j,k} dy \to 
\int_{V_0} h_j \sum_k  {{\delta_{j,k} - \varphi_j \varphi_k} \over 
{(1 + |\nabla \varphi|^2)^{3/2}}} \varphi_{j,k} dy$$
as $\epsilon \to 0,$ and  consequently  (\ref{eq(Nab1)}) also holds with $\varphi \in W^{2,2}.$ We may assume $|\nabla 
\varphi| \leq 1/2$ because $\varphi$ locally parametrizes a $C^1$ surface,
and then  (\ref{eq(Nab1)}) becomes an elliptic equation with $C^{\alpha}$ coefficients
(which depend on $\varphi$).  It then follows by Schauder's theorem (see \cite{Caff})
that $\varphi \in C^{2 + \alpha}(V_0)$ for some $\alpha > 0.$   
Finally, by the analyticity of the right side of (\ref{eq(Nab1)}),   
 the function $\varphi$, and hence the 
surface $S$, is real analytic by a theorem of
Hopf \cite{Hopf} (see also \cite{Morrey}).

\endproof

See Theorem 5 below for a related result for the case $q =1.$




 
\subsection{Radial Functions}

Assume $K$ is radial, $K(x)= K(|x|)$  
and assume $f$ is radial and $f \notin \cM_{p,q,\lambda}(f).$ 
Then    
averaging over rotations shows that every  $u \in \cM_{p,q,\lambda}(f)$ 
is radial and  
$$Du = \rho(|x|) {x \over |x|} \mu$$
where $\mu$ is invariant under rotations and where
 $\rho(|x|) = \pm 1$ a.e. $d\mu$.   
Let $H \in L^1(\mu)$ be radial and satisfy
 $\int H d\mu =0$ and $H = 0$ on $|x| < \epsilon$,  and define  
$$h(x) =  \int_{B(0,|x|)} H(|y|) {{1} 
\over {|y|^{d-1}}} d\mu.$$  
Then $h \in BV$ is radial and 
$$Dh = \vec \nu =  
 H(|x|) {x \over |x|} \mu.$$
Consequently $\vec \nu_s =0$ and 
(\ref{eq(3.1)}) gives
\begin{eqnarray*}
\int \rho H d\mu &=& \lambda \int K* J_{p,q}(x) \int_{B(0,|x|)} {{H(y)}
\over {|y|^{d-1}}} d\mu(y) dx \\
&=& \lambda \int \Bigl(\int_{|x| > |y|} K* J_{p,q}(x) dx\Bigr)
 {{H(|y|)} \over {|y|^{d-1}}} d\mu(y),
\end{eqnarray*}
 so that a.e. $d\mu$,
\begin{equation}
\rho (|y|) = {{\lambda} \over {|y|^{d-1}}} \int_{|x| > |y|}K *J_{p,q}(x)dx.
\label{eq(4.1)}
\end{equation}
 But the right side of (\ref{eq(4.1)}) is  real analytic in $|y|$, with a possible 
pole at $|y| =0,$  
 and $\rho(|y|)  = \pm 1$
almost everywhere $\mu.$  
Therefore there is a finite set 
\begin{equation}
\{r_1 < r_2 < \dots < r_n\} 
\label{eq(4.2)}
\end{equation}
of radii  such that 

$$Du = {{x} \over {|x|}} \sum_{j=1}^n c_j \Lambda_{d-1}|\{|x| = r_j\}|
$$
for real constants $c_1, \dots, c_n.$  
By Lemma \ref{Lemma2.1}, $J_{p,q}$ is uniquely determined by $f$, and hence the set (\ref{eq(4.2)}) is 
also unique.  Moreover, it
follows from Lemma \ref{Lemma2.1} that  
 for each $j$,
either $c_j \geq 0$ for all $u \in \cM_{p,1,\lambda}(f)$ or
 $c_j \leq 0$  for all $u \in \cM_{p,1,\lambda}(f)$. 
We have proved:

\begin{theorem}
\label{Thm4.1}
Suppose $K$ and $f$ are both radial. If
$f \notin \cM_{p,q,\lambda}(f)$, then there is a finite set 
(\ref{eq(4.2)}) such that all $u \in \cM_{p,q,\lambda}(f)$ 
have the form
\begin{equation}
\label{eq(4.3)}
\sum_{j=1}^n c_j \newchi_{B(0,r_j)}. 
\end{equation}
 Moreover, there is 
$X^+ \subset \{1, 2, \dots, n\}$ such that 
$c_j \geq 0$ if $j \in X^+$ while $c_j \leq 0$ if $j \notin X^+.$ 
\end{theorem}

Note that by convexity $\cM_{p,q,\lambda}(f)$ consists of a single 
function unless $p = q = 1.$ In Section \ref{RadialMinimizers} we will say more about the solutions of the form (\ref{eq(4.3)}).

\subsection{Example}

Unfortunately, Theorem \ref{Thm4.1} does not hold more generally.   The reason is 
that  when $u$ is not radial it is difficult to produce $BV$ functions satisfying 
$Dh = \vec \nu << \mu.$  For simplicity we take $d =2$ and $p = q =1$ and 
define
$$J(x,y)=\left\{
\begin{array}{rcr}
1& \mbox{ if }& 0 < x\leq 1\\ 
-1 & \mbox{ if }& -1 < x\leq 0
\end{array}
\right.
$$
\medskip
\noindent and $$J(x+2,y) = J(x,y).$$

Choose $\lambda >0$ so that $U = \lambda K *J$ satisfies $\|U\|_* =1,$
and note that
${{U} \over {|U|}} =J.$
Also notice that $u \in C^2$ solves the curvature equation
\begin{equation}{\rm {div}} \Bigl({{\nabla u} \over {|\nabla u|}}\Bigr) = U
\label{eq(5.1)}
\end{equation}
if and only if the level sets $\{u =a\}$
 are curves $y = y(x)$ that
satisfy  the simple ODE $y'' = U(x,0) (1 + (y')^2)^{3/2}$ on the line.
Consequently  (\ref{eq(5.1)}) has infinitely many solutions $u$  and 
both $u$ and $J$ satisfy (\ref{eq(3.5)}) and (\ref{eq(3.6)}).
Hence   
by Lemma \ref{Lemma2.3}, $u$ is a minimizer for $f$ provided that 
\begin{equation}
J = {{K *(f-u)} \over {|K *(f-u)|}}, 
\label{eq(5.2)}
\end{equation}
and there are many $f$ that satisfy (\ref{eq(5.2)}).  For example, one can choose $u$ and $f$ so that $f-u=J$. Note that in this example
$u$ can be real analytic except on $U^{-1}(0)$ and not piecewise constant.  Similar examples can be made 
when $(p,q) \neq (1,1).$

\section{Further Properties of Minimizers when q =1}
When $q =1$ the minimizers $u \in \cM_{p,1,\lambda}(f)$ have several additional  properties.  
The results of the next two sections do not depend on  the real analyticity of
the kernel $K$.  They also hold when $K =I$, i.e. when $\cF_{p,q,\lambda}(h) = \lambda \|h\|_p$, and 
in the case  $K =I$  somewhat stronger
results have already been proved by Allard in  \cite{Allard1}.  
However, since the arguments  in \cite{Allard1} do not apply to the case $K \neq I$ 
 we include complete but brief  proofs. 
\subsection{Layer Cake Decomposition}

Here we have been inspired by  the paper of Strang \cite{Strang}.

\begin{lemma} 
\label{Lemma6.1}
If $q =1$ and $u \in \cM_{p,1,\lambda}(f)$, then $u \in \cM_{p,1,\lambda}(u).$
\end{lemma}

\noindent {\bf {Proof:}} If 
$$\|h\|_{\dot{BV}} + \lambda\|K*(u-h)\|_p < \|u\|_{\dot{BV}},$$
then by the triangle inequality  
$$\|h\|_{\dot{BV}} + \lambda\|K*(f -h)\|_p <  \|u\|_{\dot{BV}} + \lambda\|K*(f-u)\|_p$$
so that $u$ is not a minimizer for $f$.
\endproof

We write $$\cM = \cM_{p,1,\lambda} = \bigcup_f \cM_{p,1,\lambda}(f).$$ 

\begin{lemma}
\label{Lemma6.2}
Let $u \in BV$.  Then  $u \in \cM$ if and only if 
\begin{equation}
\left| \int \rho \cdot {{d \vec \nu} \over {d \mu}} d\mu \right|
\leq \|\vec \nu_s\| + \lambda \|K * h\|_p 
\label{eq(6.1)}
\end{equation}
for all $h \in BV$, where $Dh = \vec \nu$  and $\vec \nu_s$ is 
the part of $\vec \nu$ singular to $\mu.$  
\end{lemma}

\noindent {\bf {Proof:}} By Lemma \ref{Lemma6.1} we may take $f =u.$  Then  
 for $|\epsilon|$ small we have
\begin{eqnarray*}
0 &\leq& \|u + \epsilon h\|_{\dot{BV}} - \|u\|_{\dot{BV}} + \lambda\|\epsilon K * h\|_p \\
&=& |\epsilon| \|\vec \nu_s\| + \epsilon \int \rho \cdot {{d \vec \nu} \over {d \mu}} d\mu 
+ |\epsilon| \lambda \|K * h\|_p + o(|\epsilon|) \\
\end{eqnarray*} 
and the Lemma follows from the proof of Lemma \ref{Lemma2.2}.  
\endproof

Let $a < b$ be such that 
\begin{equation}
\mu(\{u = a\} \cup \{u =b\}) =0. 
\label{eq(6.2)}
\end{equation}
 Then $u_{a,b} = {\rm {Min}}\{(u-a)^+, (b -a)\} \in BV$
and $D(u_{a,b}) = \newchi_{a < u < b} \vec {\rho}\mu.$

\begin{lemma}
\label{Lemma6.3} Assume $q =1$.

(a) If $u \in \cM$, then $u_{a,b} \in \cM.$ 

(b) More generally, if $u\in \cM$ and if $v\in BV$ satisfies  $\mu_v<<\mu_u$ and $\rho_v=\rho_u$ a.e. $d\mu_v$, then $v\in \cM$. 
\end{lemma}

\noindent {\bf Proof:} To prove (a) we verify (\ref{eq(6.1)}).  Write $\mu_{a,b} = \newchi_{(a,b)} \mu$ so that 
$D(u_{a,b}) = \vec{\rho} \mu_{a,b}.$ Let $h \in BV$ and write $Dh = \vec{\nu}.$ Then
by (\ref{eq(6.2)})
$$\vec {\nu} = \newchi_{a < u < b} {{d \vec{\nu}} \over {d \mu}} \mu  + \bigl((\vec {\nu})_s 
+ \newchi_{u(x) \notin [a,b]} {{d \vec {\nu}} 
\over {d \mu}} \mu\bigr)$$
is the Lebesgue decomposition of $\vec \nu$ with respect to $\mu_{a,b},$ and 
$$\int \vec {\rho} \cdot {{d \vec {\nu}} \over {d\mu_{a,b}}} d\mu_{a,b} =
\int \vec {\rho} \cdot {{d \vec {\nu}} \over {d\mu}} d\mu - \int_{g(x) \notin [a,b]} 
\vec {\rho} \cdot {{d \vec {\nu}} \over {d\mu}} d\mu.$$
Then (\ref{eq(6.1)}) for $\nu$ and $\mu_{a,b}$ follows from (\ref{eq(6.1)}) for $\mu$ and $\nu$.  The proof of (b) is similar. 
\endproof

For simplicity we assume $u \geq 0.$  Write $E_t =\{x: u(x) > t\}.$ Then by Evans-Gariepy \cite{EvansGariepy}, $E_t$
has finite perimeter for almost every $t$, 
\begin{equation}
\|u\|_{\dot{BV}} = \int_0^{\infty} \|\newchi_{E_t}\|_{\dot{BV}} dt, 
\label{eq(6.3)}
\end{equation}
and
\begin{equation}
u(x) = \int_0^{\infty} \newchi_{E_t}(x) dt. 
\label{eq(6.4)}
\end{equation}
Moreover, almost every set $E_t$ has a {\it {measure theoretic boundary}} $\partial_* E_t$ such 
that 
\begin{equation}
\Lambda_{d-1}(\partial_* E_t) = \|\newchi_{E_t}\|_{\dot{BV}} 
\label{eq(6.5)}
\end{equation}
and a {\it {measure theoretic outer normal}} $\vec{n}_t: \partial_* E_t \to S^{d-1}$ so that
\begin{equation}
D(\newchi_{E_t}) = \vec {n}_t\Lambda_{d-1} \bigl| \partial_* E_t.
\label{eq(6.6)}
\end{equation}


\begin{theorem}
\label{Thm6.4}
Assume $q =1$.  

(a) If $u \in \cM$, then for
 almost every $t$,
$\newchi_{E_t} \in \cM. $

(b) If $u \in \cM$ and $u \geq 0$, then for all nonnegative $c_1,...,c_n$ and for almost all $t_1<...<t_n$, $\sum c_j\chi_{E_{t_j}}\in \cM$. 
\end{theorem}

\noindent {\bf Proof:}  Suppose (a) is false.  Then there is $\beta < 1$,
and a compact set $A \subset (0,\infty)$ with  $|A| > 0$ such that for all $t \in A$
(\ref{eq(6.5)}) and (\ref{eq(6.6)}) hold and there exists
$h_t \in BV$ such that
\begin{equation}
\label{eq(6.7)}
\|\newchi_{E_t} - h_t\|_{\dot{BV}} + \lambda \|K *h_t\|_p \leq \beta \|\newchi_{E_t}\|_{\dot{BV}}. 
\end{equation}
Choose an interval $I = (a,b)$ such that (\ref{eq(6.2)}) holds and $|I \cap A| \geq {{|I|} \over {2}}.$
Define $h_t =0$ for $t \in I \setminus A,$  and take finite sums such that
\begin{equation}
\label{eq(6.8)}
\sum_{j=1}^{N_n} \newchi_{E_{t^{(n)}_j}} \Delta t_j^{(n)} \to u_{a,b} ~~ (n \to \infty), 
\end{equation}
\begin{equation}
\label{eq(6.9)}
\sum_{j=1}^{N_n} \|\newchi_{E_{t^{(n)}_j}}\|_{\dot{BV}} \Delta t_j^{(n)} \to \|u_{a,b}\| ~~ (n \to \infty), 
\end{equation}
and $t_j^{(n)} \in A$ whenever possible.  Write $h^{(n)} = \sum_{j=1}^{N_n} h_{t_j^{(n)}} \Delta t_j^{(n)}$. 
Then by (\ref{eq(6.4)}) and (\ref{eq(6.7)}) $\{h^{(n)}\}$ has a weak-star limit $h \in BV$, and by (\ref{eq(6.7)}), (\ref{eq(6.8)}) and (\ref{eq(6.9)}),
$$\|u_{a,b} -h\|_{\dot{BV}} + \lambda \|K*h\|_p \leq {{1 + \beta} \over {2}} \|u_{a,b}\|_{\dot{BV}},$$
contradicting Lemma \ref{Lemma6.3}. The proof of (b) is similar. 
\endproof

We believe that the converse of Theorem \ref{Thm6.4} is false, but we have no counterexample.   In the case $K =I$ and $p =1$ the converse of this Theorem is 
true.  See \cite{Allard1}, Theorem 5.3. 

\subsection{Characteristic Functions}

Still assuming $q =1$ we let $E$  be such that $\newchi_E \in \cM.$ Then
by Evans-Gariepy \cite{EvansGariepy} $\partial_* E = N \cup \bigcup K_j$, where $D(\newchi_E)(N) = \Lambda_{n-1}(N) =0,$
$K_j$ is compact and $K_j \subset S_j$, where $S_j$ is a $C^1-$hypersurface with continuous
unit normal $\vec{n}_j(x), x \in S_j,$ and $\vec{n}_j$ is the measure theoretic outer normal
of $E$. After a coordinate change write $S_j = \{x_d = \varphi_j(y)\}, y = (x_1, \dots, x_{d-1})$
with $\nabla \varphi_j$ continuous and ${\vec{n}_j}(y,\varphi_j(y)) =
 (1 + |\nabla \varphi_j|^2)^{-1/2} (\nabla \varphi_j, -1).$ 
Assume $y =0$ is a point of Lebesgue density of $(y,\varphi_j)^{-1}(K_j)$, let
$V \subset \R^{d-1}$ be a neighborhood of $y =0$, let $g \in C_0^{\infty}(V)$
with $g \geq 0,$  and consider the variation
$u_{\epsilon} = \newchi_{E_{\epsilon}}$
where $\epsilon > 0$ and 
$$E_{\epsilon} =  E \cup \{0 \leq x_d \leq \epsilon g(y),  y \in V\}.$$
Then $E \subset E_{\epsilon}$, and writing $u_0 = \newchi_E$, we have
\begin{equation}
\|u_{\epsilon}\|_{\dot{BV}} - \|u_0\|_{\dot{BV}} = \int_V \sqrt{(1 + |\nabla(\varphi_j + \epsilon g)|^2)}
- \sqrt{(1 + |\nabla f\varphi_j|^2)} dy + o(\epsilon)  
\end{equation}
because by \cite{EvansGariepy}  page 203 
$$\Lambda_{d-1}((\partial_* E) \cup (E_\epsilon \setminus E)) = o(\epsilon)$$  
$\Lambda_{d-1}$ a.e. on $K_j.$   
Hence
\begin{equation}
\|u_{\epsilon}\|_{\dot{BV}} - \|u_0\|_{\dot{BV}} =
 \epsilon \int_V \nabla g \cdot
 {{\nabla \varphi_j} \over {\sqrt{1 + |\nabla \varphi_j|^2}}} dy + o(\epsilon).
\label{eq(7.1)}
\end{equation}
Also,  a careful calculation gives
\begin{equation}
\lambda\|K * (u_{\epsilon} - u_0)\|_p = \lambda |\epsilon|\left\| \int_V
K(x - (y,\varphi_j(y)) g(y) dy\right\|_{L^p(dx)} + o(\epsilon). 
\label{eq(7.2)}
\end{equation}
Together (\ref{eq(7.1)}) and (\ref{eq(7.2)}) show
\begin{equation} - \int_V \nabla g \cdot \Bigl({{\nabla \varphi_j} \over {\sqrt{1 + 
|\nabla \varphi_j|^2}}}\Bigr) dy 
\leq  \lambda \|K\|_p \int_V gdy. 
\label{eq(7.3)}
\end{equation}
Repeating this argument with $\epsilon < 0$ and with $g \leq 0$  we obtain:
\begin{theorem}
\label{Thm7.1}
On the hypersurface $S_j \subset \partial_* E$
\begin{equation}
\Bigl|{\rm {div}}\Bigl({{\nabla \varphi_j} \over {\sqrt{1 + |\nabla \varphi_j|^2}}}\Bigr)\Bigr| \leq \lambda\|K\|_p. 
\label{eq(7.4)}
\end{equation}
when viewed as a distribution on $(y,\varphi_j)^{-1}(S_j).$ 
\end{theorem}
\noindent By (\ref{eq(7.4)}) and Section 7.7 of \cite{AmbrosioFuscoPallara}  
we see  that $\varphi_j \in W^{2,2}_{\rm {loc}} \cap C^{1 + \alpha}$ for any $\alpha < 1.$   
Combining Theorem \ref{Thm7.1} with Theorem \ref{Thm6.4} we obtain:
\begin{theorem}
\label{Thm7.2}
Assume $q =1$ and $u \in \cM$.  Then for almost every $t$, $E_t = \{u > t\}$ has finite
perimeter and $\Lambda_{d-1}$ almost every 
 point of the measure theoretic boundary $\partial_*{E_t}$ lies
 on a $C^{1+ \alpha}$, $\alpha < 1,$  surface having distributional mean curvature at most 
$\lambda \|K\|_p.$ 
\end{theorem}
We note that the ``distributional mean curvature"  defined by (\ref{eq(7.4)})
is the same as the generalized mean curvature defined by Allard in  \cite{Allard1},
and thus Theorem \ref{Thm7.2} complements Theorem 1.2 and Theorem 1.6 of \cite{Allard1}.
However, unlike the situation in Theorem 1, we cannot conclude that the $C^{1 + \alpha}$
surface meeting $\partial_*{E_t}$ is real analytic because  the left
side of (\ref{eq(7.4)}) may not be  H\"older continuous.









\subsection{Radial Minimizers} 
\label{RadialMinimizers}

In this section we assume $q =1$ and $p =1$.
For convenience we assume the kernel $K(x) = e^{-\pi |x|^2}$, so that $K_t$ has the form
\begin{equation}
\label{EQ7.1}
K_t(x) = t^{-d/2} K\left({{x} \over{\sqrt{t}}}\right) 
\end{equation}
and 
\begin{equation}
\label{EQ7.2}
K_s * K_t = K_{s+t}. 
  \end{equation}             
Note that 
(\ref{EQ7.1}) and (\ref{EQ7.2}) imply that  
\begin{equation}
\label{EQ7.3}
\|K_t * f\|_1 ~~~decreases~ in ~ t 
\end{equation}
and for $f \in L^1$ with compact support 
\begin{equation}
\label{EQ7.4}
\lim_{t \to \infty}\|K_t * f\|_1 = \left|\int f dx\right|.
\end{equation}
For fixed $\lambda$ and $t$ we set $$R(\lambda,t) = \{r > 0:
\newchi_{B(0,r)} \in \cM\}.$$ By Theorem \ref{Thm4.1} and Theorem \ref{Thm6.4}, we have $R(\lambda,t) \neq \emptyset.$
For $t =0$  and  $K =I$  our problem (\ref{(1.2)}) becomes the problem 
$$\inf \{\|u\|_{\dot{BV}} + \lambda \|f - u\|_{L^1}\}$$ studied by  Chan and Esedoglu  in \cite{ChanEsedoglu},
and in that case Chan and Esedoglu showed $R(\lambda,0) = [{{2} \over {\lambda}},
\infty).$

\begin{theorem}
\label{THM7.1}
There exists $r_0 = r_0(\lambda,t)$ such that
\begin{equation}
\label{EQ7.5}
R(\lambda,t) = [r_0,\infty).
\end{equation}
Moreover  
\begin{equation}
\label{EQ7.6}
[0,\infty) \ni t \to r_0(t) ~~is~~ nondecreasing 
\end{equation}
and
\begin{equation}
\label{EQ7.7}
\lim_{t \to \infty} r_0(t) = \infty. 
\end{equation}
\end{theorem}

\noindent {\bf Proof:} Assume $r \notin R(\lambda,t)$
and $0 < s < r.$ Write $\alpha = {{r} \over {s}} > 1$ and $f = \newchi_{B(0,r)}.$ By hypothesis there is $g \in BV$ such that
\begin{equation}
\label{EQ7.8}
\|g\|_{\dot{BV}} + \lambda \|K_t*(f-g)\|_1 < \|f\|_{\dot{BV}}.
\end{equation}
We write $\tilde g(x) = g(\alpha x), \tilde f(x) = f(\alpha x) = \newchi_{B(0,s)}(x),$ and 
change variables carefully in (\ref{EQ7.8}) to get
$$\alpha \|\tilde g\|_{\dot{BV}} + \lambda\|{{1} \over {t^{d/2}}} \int K\bigl({{x -y} \over {\sqrt{t}}}\bigr)
(\tilde f - \tilde g)\bigl({{y} \over {\alpha}}\bigr) dy\|_{L^1(x)} < \alpha \|\tilde f\|_{\dot{BV}}$$
so that
$$\alpha \|\tilde g\|_{\dot{BV}} + \lambda \|{{\alpha^d} \over {t^{d/2}}} \int K\bigl({{\alpha x' - \alpha y'} 
\over {\sqrt{t}}}\bigr) (\tilde f - \tilde g)(y') dy'\|_{L^1(\alpha x')} < \alpha \|\tilde f\|_{\dot{BV}}$$
and 
$$\alpha \|\tilde g\|_{\dot{BV}} + 
\lambda \alpha^d \int\Bigl|K_{{t} \over {\alpha^2}} * (\tilde f - \tilde g)(x')\Bigr| dx' 
< \alpha \|\tilde f\|_{\dot{BV}}.$$
Since $\alpha > 1$, this and (\ref{EQ7.3}) show 
$$\|\tilde g\|_{\dot{BV}} + \lambda\| K_{t}* (\tilde f - \tilde g)\|_1
< \|\tilde f\|_{\dot{BV}}
$$
so that $s \notin R(\lambda,t).$ That proves (\ref{EQ7.5}), and (\ref{EQ7.6}) now follows easily
from (\ref{EQ7.3}). To prove (\ref{EQ7.7}) take $g = {{r^d} \over {s^d}} \newchi_B(0,s),~ s >r$
and use (\ref{EQ7.4}).
\endproof


We note that not all radial minimizers have the form $\newchi_{B(0,r)}$. This is seen by considering separately, for large fixed $t$ and $\lambda$,  the 
function  
$\newchi_{B(0,r_2)} + \newchi_{B(0,r_1)}
$ with $r_1$ and $r_2 - r_1$ large.   

\end{document}